\documentclass[12pt]{amsart}
\usepackage{amssymb,amsmath,amsfonts,latexsym,amscd}
\usepackage{bm}
\usepackage{ifthen}
\usepackage{graphicx}
\usepackage{color}
\newcommand{\vp}{\varphi}

\newcommand{\clb}{\mathcal{B}}

\newcommand{\clh}{\mathcal{H}}

\newcommand{\clm}{\mathcal{M}}

\newcommand{\raro}{\rightarrow}

\newcommand{\be}{\begin{equation}}
\newcommand{\ee}{\end{equation}}
\newcommand{\ben}{\begin{eqnarray*}}
\newcommand{\een}{\end{eqnarray*}}

\newcommand{\bi}{\begin{itemize}}
\newcommand{\ei}{\end{itemize}}
\newcommand{\es}{{\mathcal S}}

\newcommand{\D}{\mathbb{D}}

\newcommand{\C}{\mathbb{C}}

\newcommand{\clq}{\mathcal{Q}}
\newcommand{\cls}{\mathcal{S}}

\newcommand{\ds}{\displaystyle}

\newtheorem{thm}{Theorem}[section]
\newtheorem{cor}[thm]{Corollary}
\newtheorem{lem}[thm]{Lemma}
\newtheorem{prop}[thm]{Proposition}
\newtheorem{rem}[thm]{Remark}

\newtheorem{eg}[thm]{Example}
\newtheorem{qn}[thm]{Question}
\newenvironment{pf}[1][]{%
 \vskip 3mm
 \noindent
 \ifthenelse{\equal{#1}{}}%
  {{\slshape Proof. }}%
  {{\slshape #1.} }%
 }%
{\qed\bigskip}

\newcommand{\blem}{\begin{lem}}
\newcommand{\elem}{\end{lem}}
\newcommand{\bthm}{\begin{thm}}
\newcommand{\ethm}{\end{thm}}
\newcommand{\bprop}{\begin{prop}}
\newcommand{\eprop}{\end{prop}}
\newcommand{\bpf}{\begin{pf}}
\newcommand{\epf}{\end{pf}}
\newcommand{\bcor}{\begin{cor}}
\newcommand{\ecor}{\end{cor}}
\newcommand{\beg}{\begin{eg}}
\newcommand{\eeg}{\end{eg}}
\newcommand{\brem}{\begin{rem}}
\newcommand{\erem}{\end{rem}}
\newcommand{\bei}{\begin{itemize}}
\newcommand{\eei}{\end{itemize}}
\newcommand{\bqn}{\begin{qn}}
\newcommand{\eqn}{\end{qn}}

\begin{document}
\title[Model spaces invariant under composition operators]
{Model spaces invariant under composition operators}

\author[Muthukumar]{P. Muthukumar}
\address{P. Muthukumar, Indian Statistical Institute,
Statistics and Mathematics Unit, 8th Mile, Mysore Road,
Bangalore, 560 059, India.}
\email{pmuthumaths@gmail.com}

\author[Sarkar]{Jaydeb Sarkar}
\address{Jaydeb Sarkar, Indian Statistical Institute,
Statistics and Mathematics Unit, 8th Mile, Mysore Road,
Bangalore, 560 059, India.}
\email{jaydeb@gmail.com, jay@isibang.ac.in}



\begin{abstract}
Given a holomorphic self-map $\varphi$ of $\D$ (the open unit disc in $\mathbb{C}$), the composition operator $C_{\varphi} f = f \circ \varphi$, $f \in H^2(\mathbb{\D})$, defines a bounded linear operator on the Hardy space $H^2(\mathbb{\D})$. The model spaces are the backward shift-invariant closed subspaces of $H^2(\mathbb{\D})$, which are canonically associated with inner functions. In this paper, we study model spaces that are invariant under composition operators. Emphasis is put on finite-dimensional model spaces, affine transformations, and linear fractional transformations.
\end{abstract}

\maketitle

\section{Preliminaries}\label{Mjay2-sec2}

Let $\D$ denote the open unit disk in $\mathbb{C}$ and let $H^2(\D)$ be the Hardy space over $\mathbb{D}$. For $\vp$ a holomorphic self-map of $\D$, define the \textit{composition operator} $C_\vp$ on $H^2(\D)$ by
\[
C_\vp f = f \circ \vp \qquad (f \in H^2(\D)).
\]
It follows from the Littlewood subordination principle that $C_\vp \in \clb(H^2(\D))$ \cite{Shapiro:Book}. Here $\clb(H^2(\D))$ denote the set of all bounded linear operators on $H^2(\D)$.

The theory of composition operators is highly interdisciplinary with its natural connections to complex analysis, linear dynamics, complex geometry, and functional analysis. For instance, the question of the existence of nontrivial invariant subspaces of bounded linear operators on Hilbert spaces can be formulated in terms of minimal invariant subspaces of composition operators induced by hyperbolic automorphisms of $\D$ \cite{MPS, ISP-CO1, ISP-CO}.

The theme of this paper is also in line with the invariant subspaces of composition operators, but restricted to model (or inner function based quotient) spaces of $H^2(\D)$. Recall that a function $\theta \in H^\infty(\D)$ is said to be \textit{inner} if $|\theta| = 1$ a.e. on $\partial \D$ in the sense of nontangential boundary values, where $H^\infty(\D)$ denotes the commutative Banach algebra of all bounded holomorphic functions on $\D$. Given an inner function $\theta$, the \textit{model space} $\clq_\theta$ is the quotient space
\[
\clq_\theta = H^2(\D) \ominus \theta H^2(\D) \cong H^2(\D)/ \theta H^2(\D).
\]
The aim of this paper is to describe some results and methods of the calculation of model spaces that are invariant under composition operators. Note, on the other hand, a celebrated result of A. Beurling \cite{InvSub-MO} states that a closed subspace $\clq \subseteq H^2(\D)$ is invariant under $M_z^*$ if and only if $\clq$ is a model space. Here $M_zf = z f$, $f \in H^2(\D)$, and $M_z^*$ is the backward shift
\[
M_z^* f= \frac{f(z) - f(0)}{z} \qquad (f \in H^2(\D)).
\]

In other words, we are interested in joint invariant subspaces of the (non commuting) pair $(M_z^*, C_\vp)$ on $H^2(\D)$. Yet another motivation comes from the recent classification \cite{ISCO} of joint invariant subspaces of the non commuting pair $(M_z, C_\vp)$ on $H^2(\D)$. In this context, we further remark that the concept of joint invariant subspaces of $(M_z, C_\vp)$ was independently introduced by Chalendar and Partington \cite{CP} and Mahvidi \cite{Inv2001}. For recent results related to the subject we refer to \cite{ISCO, Inv2005, Inv2015} and the references therein.

Often results on composition operators are case-based. For instance, the theory of $C_{\vp}$ corresponding to holomorphic self-map of $\D$ as
\[
\vp(z) = az+b \text{ or } \frac{az+b}{cz+d},
\]
is extremely demanding and complete answers to many basic questions are not known (cf. \cite{Deddens, Adjoint}, and also see \cite{Eva, Pascal} for more recent accounts). Even norm computations and invariant subspaces of $C_\vp$ for linear fractional transformations are known to be case-based (cf. \cite{Bourdon1, MPS}). In fact, the results of this paper are also no exception. For example, Theorem \ref{affine} states: Let $\vp$ be a holomorphic self-map of $\D$. Fix $\alpha \in \D$ and $n \in \mathbb{N}$, and consider the inner function
\[
\theta(z)=\left(\frac{z- \alpha}{1-\bar{\alpha}z}\right)^n \qquad (z \in \D).
\]
Then $\clq_\theta$ is invariant under $C_\vp$ if and only if there exist scalars $a,b$, and $c (\neq 0)$ such that
\begin{equation}\label{eqn: sect 1}
\vp(z)= \begin{cases}
a + bz & \mbox{if } \alpha = 0 \text{ and } n > 1 \\
\frac{1-c}{\bar{\alpha}} + cz & \mbox{if } \alpha \neq 0,
\end{cases}
\end{equation}
where, $|a| + |b| \leq 1$ and $|\frac{1-c}{\bar{\alpha}}| + |c| \leq 1$.

Note that the above model space $\clq_\theta$ is an $n$-dimensional Hilbert space. Even though we all know the structure of finite-dimensional model spaces (they are essentially parameterized by $n$-Blaschke products) but we do not know a general method to classify holomorphic self-maps $\vp$ of $\D$ for which a given finite-dimensional model space remains invariant under $C_\vp$ (see Question \ref{Quest: 1}). In order to make this problem more convincing, in Example \ref{example:1} we consider the $2$-dimensional model space $\clq_\theta$ corresponding to the (special) Blaschke product $\theta(z) = z \frac{z- \alpha}{1 - \bar{\alpha} z}$, $z \in \D$, and prove that $C_{\vp} \clq_\theta \subseteq \clq_\theta$ if and only if $\vp$ is the special M\"{o}bius map
\[
\vp(z) = \frac{(c_1 - 1) + (\bar{\alpha} + c_2) z}{\bar{\alpha}(c_1 + c_2z)} \qquad (z \in \D),
\]
for some scalars $c_1$ and $c_2$ such that at least one of them is nonzero (see \eqref{eqn: phi = Mob}). The above answers do not seem to tell us much about a possible answer to our general question.

We also put special emphasis on composition operators induces by linear fractional symbols (cf. Theorem \ref{modelinv}):
Let $\theta$ be an inner function, and let $\vp(z)= \frac{az+b}{cz+d}$ be a nonconstant holomorphic self-map of $\D$, where $ad - bc = 1$. Suppose $\sigma(z)= \frac{\bar{a}z-\bar{c}}{-\bar{b}z+\bar{d}}$ and $\theta(0) = 0$. If we denote $\psi = (\frac{\theta}{z}\circ\sigma)(\frac{\theta}{z})^{-1}$, then the following are equivalent:
\begin{enumerate}
  \item $\clq_\theta$ is invariant under $C_\vp$.
  \item $\psi \in H^2(\D)$ or $\psi \in H^\infty(\D)$.
  \item $\frac{\theta}{z}H^2(\D)$ is invariant under $C_\sigma$.
\end{enumerate}
A similar (but not same as above) consideration also holds for composition operators induced by affine transformations (see Theorem \ref{tm: FLT 1}).

In the application part, we consider direct interpretations of our results to model spaces and reducing subspaces of composition operators. For instance, Theorem \ref{thm: all quotient inv sub} states that:
\begin{enumerate}
\item $\clq_\theta$ is invariant under $C_\vp$ for all inner function $\theta$ if and only if $\vp$ is the identity map.
\item $\clq_\theta$ is invariant under $C_\vp$ for all holomorphic self-map $\vp$ if and only if $\theta (z)=\alpha z$ or $\theta\equiv\alpha$, where $\alpha$ is an unimodular constant.
\end{enumerate}

We now discuss the issue of reducing subspaces of composition operators from our perspective. Let $T$ be a bounded linear operator on a Hilbert space $\clh$. A closed subspace $\mathcal{M} \subseteq \clh$ is called reducing subspace for $T$ (or $\clm$ reduces $T$) if both $\mathcal{M}$ and $\mathcal{M}^{\bot}$ are invariant under $T$. It is known that $M_z$ on $H^2(\D)$ is irreducible, that is, $M_z$ has no nontrivial reducing subspaces. In the setting of composition operators, if $\alpha \in \C$ with $\text{Re}(\alpha) > 0$, then $C_{\vp_{\alpha}}$ is also irreducible \cite[Theorem 1.3]{MPS}, where
\[
\vp_{\alpha}(z)=\frac{(2- \alpha )z+ \alpha}{-\alpha z+(2+ \alpha)} \qquad (z \in \D).
\]
Another example is the composition operator corresponding to univalent loxodromic type self-map with nonzero Denjoy-Wolff point \cite[Theorem 1]{Guyker}. At the other extreme, we discuss reducibility of compositor operators restricted to model spaces. In Theorem \ref{thm: reducing}, we prove the following: Let $\vp$ be a holomorphic self-map of $\D$, $\alpha \in \D$, $n \geq 1$, and suppose $\theta(z)=\left(\frac{z- \alpha}{1-\bar{\alpha}z}\right)^n$. Then $\clq_\theta$ reduces $C_\vp$ if and only if there exist a scalar $c$, $|c| \leq 1$, and $\psi \in H^\infty(\D)$, $\|\psi\|_{\infty} \leq 1$, such that
\[
\vp(z)= \begin{cases}
z \psi & \mbox{if } \alpha = 0 \text{ and } n = 1 \\
c z & \mbox{if } \alpha = 0 \text{ and } n \geq 2\\
z & \mbox{if } \alpha \neq 0.
\end{cases}
\]

Finally, we make a general remark. Recall that a linear fractional transformation $\vp(z) = \frac{az+b}{cz+d}$ is a self-map of $\D$ if and only if
\[
|b \bar{d} - a \bar{c}| + |ad - bc| \leq |d|^2 - |c|^2.
\]
In particular, since $\vp$ is a holomorphic self-map, the inequalities following \eqref{eqn: sect 1} are automatically true. In this paper, whenever no confusion can arise, we will not explicitly mention the above inequality for linear fractional transformations that are self-maps of $\D$. For instance, the above inequality is applicable to Theorem \ref{affine}, Example \ref{example:1}, and Theorem \ref{tm: FLT 1}.

The rest of the paper is organized as follows. Section \ref{Sec: Beurling type} recalls (rather a modified version of) the classification of Beurling type invariant subspaces from \cite{ISCO}. Section \ref{Mjay2-sec4} relates affine transformations and certain finite-dimensional model spaces. Section \ref{sec: LFS} deals with model subspaces invariant under composition operators corresponding to fractional linear transformations. The final section, Section \ref{Mjay2-sec5}, deals with some of the direct applications of our main results. We also discuss irreducible composition operators at the level of model spaces.

\section{Beurling type invariant subspaces}\label{Sec: Beurling type}

We begin with the formal definition of Hardy space $H^2(\D)$. Recall that the Hardy space is the function Hilbert space of all holomorphic functions on $\D$ whose power series have square-summable coefficients. The following equivalent description of $H^2(\D)$ is well known:
\[
H^2(\D) = \Big\{f \in \text{Hol}(\D): \|f\|:= \Big[\sup_{0 \leq r < 1} \frac{1}{2 \pi} \int_{0}^{2\pi} |f(r e^{it})|^2 dt\Big]^{\frac{1}{2}} < \infty\Big\}.
\]
Also we denote the closed unit ball of $H^\infty(\D)$ (which is popularly known as Schur class) by
\[
\es(\D) = \{\psi \in H^\infty(\D): \|\psi\|_\infty: = \sup_{z \in \D} |\psi(z)| \leq 1\}.
\]
We follow the standard notation: Given a bounded linear operator $T$ on some Hilbert space $\clh$, the lattice of all closed $T$-invariant subspaces is denoted by $\text{Lat} T$, that is
\[
\text{Lat} T = \{\cls \subseteq \clh \text{ closed subspace}: T \cls \subseteq \cls\}.
\]
By Beurling's Theorem, we know that
\[
\text{Lat} M_z = \{\{0\}, \theta H^2(\D): \theta \text{ inner}\}, \text{ and } \text{Lat} M_z^* = \{H^2(\D), \clq_\theta: \theta \text{ inner}\}.
\]
A closed subspace $\cls \subseteq H^2(\D)$ is said to be a \textit{Beurling type invariant subspace} if $\cls = \theta H^2(\D)$ for some inner function $\theta$, or equivalently (in view of Beurling), $\cls$ is invariant under $M_z$. The purpose of this short section is to recall the classification of Beurling type invariant subspaces of composition operators  \cite[Theorem 2.3]{ISCO}. However, the present version is a simple refinement of \cite{ISCO}, and we only prove the new component of the result.

\bthm\label{main}
Let $\theta$ be an inner function and $\vp$ be a holomorphic self-map of $\mathbb{D}$. The following are equivalent:
\begin{enumerate}
 \item $\theta H^2(\D) \in \text{Lat} C_\vp$.
 \item $\ds\frac{\theta\circ\vp}{\theta}\in \es(\D)$.
 \item $\ds\frac{\theta\circ\vp}{\theta}\in H^\infty(\D)$.
 \item $\ds\frac{\theta\circ\vp}{\theta}\in H^2(\D)$.
 \end{enumerate}
\ethm
\bpf
The equivalence $(1) \Leftrightarrow (2)$ follows from \cite[Theorem 2.3]{ISCO}. Since
\[
\es(\D) \subsetneq H^\infty(\D) \subsetneq H^2(\D),
\]
it follows that $(2) \Rightarrow (3) \Rightarrow (4)$. Therefore, it suffices to verify that $(4) \Rightarrow (2)$. Suppose that $\frac{\theta\circ\vp}{\theta}\in H^2(\D)$. Since $\theta$ is inner, we have
\[
\left|\frac{\theta\circ\vp}{\theta}(e^{it})\right|= |(\theta\circ\vp)(e^{it})|\leq 1 \mbox{~for all~}
t\in [0,2 \pi] \mbox{~a.e.}
\]
Therefore, $\frac{\theta\circ\sigma}{\theta}  \in H^\infty(\D)$ with $\left\|\frac{\theta\circ\vp}{\theta}\right\|_\infty\leq 1$ \cite[Theorem 2.11]{Duren:Hpspace} , that is, $\frac{\theta\circ\vp}{\theta}\in \es(\D)$, which completes the proof of the theorem.
\epf

This result will play a key role in what follows.

\section{Finite-dimensional model spaces}\label{Mjay2-sec4}

Recall that the Hardy space $H^2(\D)$ is associated with the Szeg\"{o} kernel $k: \D \times \D \raro \mathbb{C}$, where
\[
k(z, w) = (1 - z \bar{w})^{-1} \qquad (z, w \in \D).
\]
Therefore, the linear span of $\{k(\cdot, w): w \in \D\}$ is dense in $H^2(\D)$, and $f(w) = \langle f, k(\cdot, w) \rangle_{H^2(\D)}$ for all $f \in H^2(\D)$ and $w \in \D$. For each $\alpha \in \D$, define the \textit{Blaschke factor} $b_{\alpha}$ by
\[
b_{\alpha}(z) = \frac{z - \alpha}{1 - \bar{\alpha} z} \qquad (z \in \D).
\]
It is known that a model space $\clq_{\theta}$ is finite dimensional if and only if $\theta$ is a finite Blaschke product (cf. \cite[Theorem 3.8]{Radjavi book}). More specifically, if $\text{dim}\clq_{\theta} = n < \infty$, then there exists $\alpha_1, \ldots, \alpha_m \in \D$, and $n_1, \ldots, n_m \in \mathbb{N}$ such that $\sum_{i=1}^m n_i = n$ and
\[
\theta =\prod_{i=1}^{m} b_{\alpha_i}^{n_i}.
\]
Moreover, $\{c_{\alpha_i}^{(l_i)}: 0 \leq l_i \leq n_i-1, 1 \leq i \leq m\}
$
is a basis for $\clq_{\theta}$ (see \cite[Proposition 5.16]{Ross book}),  where
\[
c_{\alpha}^{(t)} := \frac{z^{t}}{(1 - \bar{\alpha} z)^{t+1}} \qquad (t \in \mathbb{Z}_+).
\]
In the following, we prove that a ``simple'' finite dimensional models space is invariant under $C_\vp$ if and only if $\vp$ is a special affine map.

\bthm\label{affine}
Let $\vp$ be a holomorphic self-map of $\D$, $\alpha \in \D$, $n \in \mathbb{N}$, and suppose $\theta(z)=\left(\frac{z- \alpha}{1-\bar{\alpha}z}\right)^n$. Then $\clq_\theta \in \text{Lat} C_\vp$ if and only if there exist scalars $a,b$, and $c (\neq 0)$ such that
\[
\vp(z)= \begin{cases}
a + bz & \mbox{if } \alpha = 0 \text{ and } n > 1 \\
\frac{1-c}{\bar{\alpha}} + cz & \mbox{if } \alpha \neq 0.
\end{cases}
\]
\ethm
\bpf
In view of the preceding discussion, we have
\[
\clq_\theta = \text{span}\left\{ \frac{p(z)}{(1-  \bar{\alpha} z)^n}: p \in \C[z], \text{deg } p \leq n-1 \right\}.
\]
Suppose $\alpha = 0$. Then $\theta(z)=z^n$ and
\[
\clq_\theta=\mbox{span}\{1, z, z^2,\ldots,z^{n-1}\}.
\]
Let $\vp = a + bz$. Clearly, if $f = \sum_{j=0}^{n-1} a_{j} z^{j} \in \clq_\theta$, then
\[
f\circ \vp = \sum_{j=0}^{n-1} a_{j} \vp^{j} = \sum_{j=0}^{n-1} a_{j} (a + bz)^{j},
\]
and hence $C_\vp f = f\circ \vp \in \clq_\theta$. For the converse part, suppose $C_{\vp} \clq_\theta \subseteq \clq_{\theta}$. Since the identity map $f(z)=z$, $z \in \D$, is in $\clq_\theta$, if follows that
\[
f\circ \vp=\vp \in \clq_\theta,
\]
which implies that $\vp$ is a polynomial. Also note that $g \in \clq_\theta$, where $g(z)=z^{n-1}$, $z \in \D$. Then
\[
g\circ \vp=\vp^{n-1} \in \clq_\theta,
\]
which implies that
\[
(n-1) \mbox{deg}(\vp)=\mbox{deg}(\vp^{n-1})\leq n-1.
\]
Since $n-1\neq 0$, it follows that $\mbox{deg}(\vp)\leq 1$, and hence $\vp$ is an affine map. Now suppose that $\alpha \neq 0$, and let $C_{\vp} \clq_\theta \subseteq \clq_\theta$. Since $k(\cdot, \alpha) \in \clq_\theta$, it follows that
\[
k(\cdot, \alpha) \circ \vp=\ds\frac{1}{1-\bar{\alpha}\vp} \in \clq_\theta,
\]
and hence, there exists a polynomial $p$ of degree $\leq n-1$ such that
\[
\frac{1}{1- \bar{\alpha} \vp}=\frac{p(z)}{(1-\bar{\alpha}z)^n}.
\]
In particular, $\vp$ is a nonconstant polynomial. Suppose $\mbox{deg}(\vp)=m >0$. Since
\[
g = \frac{1}{(1-\bar{\alpha}z)^n} \in \clq_\theta,
\]
we have
\[
g\circ \vp = \frac{1}{(1-\bar{\alpha}\vp)^n}\in \clq_\theta.
\]
Thus, there exists $q \in \C[z]$ such that
\[
(1-\bar{\alpha}\vp)^n q(z)=(1-\bar{\alpha}z)^n.
\]
Comparing the degrees of both sides, we have $mn\leq n$, which implies $\mbox{deg}(\vp)=1$. Then $\mbox{deg}(q)=0$, and hence $q$ is a nonzero constant. Therefore, there exists a nonzero scalar $c$ such that $1-\bar{\alpha}\vp=c (1-\bar{\alpha}z)$, that is
\[
\vp(z)= \frac{1-c}{\bar{\alpha}} + cz.
\]
We now turn to the converse direction. Let $c \neq 0$ and $\vp(z)= \frac{1-c}{\bar{\alpha}}+cz$. Then $1 - \bar{\alpha} \vp = c (1 - \bar{\alpha} z)$, and hence
\[
\left(\frac{1}{1- \bar{\alpha} \vp} \right)^l\in \mbox{span}\left\{\frac{1}{(1- \bar{\alpha} z)^j}: j =1, \ldots, n \right\}
= \clq_\theta
\]
for all $l=1,2,\ldots, n$. It follows that $f\circ\vp \in \clq_\theta$ for all $f\in \clq_\theta$. This completes the proof of the theorem.
\epf

We would like to point out that even affine transformations (as symbols) play a prominent role in the complexity of composition operators (cf. Deddens \cite{Deddens}).

\brem\label{constant}
In the context of Theorem \ref{affine}, we note that if $\alpha = 0$ and $n = 1$, then $\clq_z=(z H^2(\D))^{\bot}$, which is the one dimensional subspace of all constant functions in $H^2(\D)$. Clearly, $\clq_z$ is invariant under every composition operator.
\erem

We do not know the answer to the following compelling question.

\bqn\label{Quest: 1}
Consider the finite Blaschke product $\prod_{i=1}^{m} b_{\alpha_i}$ corresponding to $\alpha_1, \ldots, \alpha_n \in \D$. Characterize holomorphic self-maps $\vp$ of $\D$ such that
\[
C_{\vp} \clq_{\prod_{i=1}^{m} b_{\alpha_i}} \subseteq \clq_{\prod_{i=1}^{m} b_{\alpha_i}}.
\]
\eqn

We also put forward the following simpler (two-dimensional) version: Suppose $\alpha \neq \beta$ in $\D$, and suppose $\theta= b_{\alpha} b_{\beta}$. Classify holomorphic self-maps $\vp$ of $\D$ such that $C_\vp \clq_\theta \subseteq \clq_\theta$.

We do not know the complete answer to this question but the following special case ($\alpha \neq 0$ and $\beta = 0$):

\beg\label{example:1}
Let $\alpha \in \D$ be a nonzero scalar, and suppose $\theta = z b_{\alpha}$. Let $\vp$ be a holomorphic self-map of $\D$. Then $C_{\vp} \clq_\theta \subseteq \clq_\theta$ if and only if $\vp$ is the special M\"{o}bius map
\begin{equation}\label{eqn: phi = Mob}
\vp(z) = \frac{(c_1 - 1) + (\bar{\alpha} + c_2) z}{\bar{\alpha}(c_1 + c_2z)} \qquad (z \in \D),
\end{equation}
for some scalars $c_1$ and $c_2$ such that at least one of them is nonzero. Indeed, here the $2$-dimensional model space $\clq_\theta$ is given by $\clq_\theta = \text{span} \{1, k(\cdot, \alpha)\}$. As the space of all constant functions is invariant under every composition operators (see Remark \ref{constant}), $\clq_\theta$ is invariant under $C_\vp$ if and only if
\[
k(\cdot, \alpha) \circ \vp = \frac{1}{1 - \bar{\alpha} \vp} \in \clq_\theta.
\]
This is equivalent to the fact that
\[
\frac{1}{1 - \bar{\alpha} \vp} = \frac{c_1 + c_2 z}{1 - \bar{\alpha} z} \qquad (z \in \D),
\]
for some scalars $c_1$ and $c_2$ such that at least one of them is nonzero, which is equivalent to \eqref{eqn: phi = Mob}. This completes the proof of the claim.
\eeg

Putting together the representations of linear fractional transformations from \eqref{eqn: phi = Mob} and Theorem \ref{affine}, it is perhaps evident that the answer to Question \ref{Quest: 1} will be also case-based.

\section{Linear fractional symbols}\label{sec: LFS}

Our aim in this section is to describe the structure of model invariant subspaces of composition operators whose symbols are linear fractional transformations mapping $\D$ into itself. We are also motivated by the intricate analysis presented in Cowen \cite{Adjoint} and Deddens \cite{Deddens}.

Before we begin with a general observation, we introduce the following (standard) notation: For each $\vp \in H^\infty(\D)$, we denote by $M_\vp$ the multiplication operator $M_\vp f = \vp f$, $f \in H^2(\D)$. It is well known that $M_\vp \in \clb(H^2(\D))$ with $\|M_\vp\| = \|\vp\|_{\infty}$ for all $\vp\in H^\infty(\D)$, and $\{M_z\}' = \{M_{\vp}: \vp \in H^\infty(\D)\}$, where  $\{M_z\}'$ denotes the commutator of $M_z$.

Let $\vp$ be a holomorphic self-map of $\D$ and let $\theta$ be an inner function. Clearly, $C_\vp \clq_\theta \subseteq \clq_\theta$ if and only if
$C_\vp^*(\theta H^2(\D)) \subseteq \theta H^2(\D)$, which is, by the Douglas range inclusion theorem \cite[Theorem 1]{Douglas:Inclusion} (as $\mbox{ran} (C_{\vp}^* M_{\theta}) \subseteq \mbox{ran} M_{\theta}$), equivalent to the existence of some $X \in \clb(H^2(\D))$ such that
\[
C_\vp^* M_\theta= M_\theta X.
\]
Now if the above equality holds, then $\theta X(f) = C_{\vp}^*(\theta f)$ for all $f \in H^2(\D)$, and hence zeros of $\theta$ are also zeros of $C_{\vp}^*(\theta f)$ and
\[
\text{multiplicity}_{w} \theta \leq \text{multiplicity}_{w} C_{\vp}^*(\theta f),
\]
for all $w \in \mathcal{Z}(\theta)$ (the zero set of $\theta$). This implies that $\frac{1}{\theta} C_\vp^*(\theta f)$ is a well defined holomorphic function on $\D$. Therefore, we have proved the following general (but perhaps less practical) characterization.

\blem\label{eqn: Xf}
$Q_{\theta}$ is invariant under $C_\vp$ if and only if $f \in H^2(\D) \mapsto \frac{1}{\theta} C_\vp^* M_\theta f$ defines a bounded linear operator on $H^2(\D)$.
\elem

Now we turn to linear fractional transformations. The following adjoint formula is due to Cowen \cite[Theorem 2]{Adjoint}. There is a wide range of applications of Cowen's adjoint formula (cf. \cite{Dechao, Eva 1, Eva}).

\bthm\label{adjoint}
Let $\vp(z)=\frac{az+b}{cz+d}$ be a nonconstant holomorphic self-map of $\D$, where $ad - bc = 1$.
Then $\sigma(z) = \frac{\bar{a}z-\bar{c}} {-\bar{b}z+\bar{d}}$ is a holomorphic self-map of $\D$, $g(z)= \frac{1}{-\bar{b}z+\bar{d}}$ and $h(z)=cz+d$ are in $H^\infty(\D)$, and $C_\vp^*=M_g C_\sigma M_h^*$.
\ethm

Note that if $\vp$ is a constant map, then $\clq_\theta$ is invariant for $C_\vp$ if and only if $\theta(0)=0$. Indeed, the constant function $1 \in \clq_\theta$ if and only if $\theta(0) = \langle 1, \theta \rangle_{H^2(\D)} = 0$. Now we turn to nonconstant affine self-maps.

In what follows, we denote by $\clh(\D)$ either $H^2(\D)$ or $H^\infty(\D)$.

\bthm\label{tm: FLT 1}
Let $\theta$ be an inner function, and let $\vp(z)=az+b, a\neq0$. Suppose $\sigma(z)= \frac{\bar{a}z}{1-\bar{b}z}$, $z \in \D$. The following are equivalent.
\begin{enumerate}
\item $\clq_\theta \in \text{Lat} C_\vp$.
\item $\frac{\theta\circ\sigma}{\theta}  \in \clh(\D)$.
\item $\theta H^2(\D) \in \text{Lat} C_\sigma$.
\end{enumerate}
\ethm
\bpf
If we apply Theorem \ref{adjoint} to the case when $c=0$ and $d = 1$, we get
\[
C_\vp^*=M_\psi C_\sigma,
\]
where $\psi(z)= \frac{1}{1-\bar{b}z}$. Now suppose that $\clq_\theta$ is invariant under $C_\vp$. By Lemma \ref{eqn: Xf}, we know that the map
\[
Xf= \frac{1}{\theta} C_\vp^*M_\theta f \qquad (f \in H^2(\D)),
\]
defines a bounded linear operator $X$ on $H^2(\D)$, and hence, using $C_\vp^*=M_\psi C_\sigma$, we have
\[
Xf = \frac{1}{\theta} C_\vp^*M_\theta f = \frac{1}{\theta} \psi (\theta\circ\sigma) (f\circ\sigma) \in H^2(\D),
\]
for all $f\in H^2(\D)$. Applying this with $f = 1$ yields
\[
\psi\frac{\theta\circ\sigma}{\theta}  \in H^2(\D).
\]
But since $\frac{1}{\psi}= 1 - \bar{b}z \in H^\infty(\D)$, it follows that $\frac{\theta\circ\sigma}{\theta} \in H^2(\D)$. Theorem \ref{main} then implies that $\theta H^2(\D)$ is invariant under $C_\sigma$. This proves $(1) \Rightarrow (2) \Rightarrow (3)$.

\noindent Now we assume that $C_\sigma (\theta H^2(\D)) \subseteq \theta H^2(\D)$. Again, Theorem \ref{main} implies that $\frac{\theta\circ\sigma}{\theta} \in \es(\D)$, and therefore
$\psi \frac{\theta\circ\sigma}{\theta}  \in H^\infty(\D)$. Then $\frac{1}{\theta} C_{\vp}^* M_\theta$ defines a bounded linear operator on $H^2(\D)$, and hence, $\clq_\theta$ is invariant under $C_\vp$. This proves that $(3) \Rightarrow (2) \Rightarrow (1)$ and completes the proof of the theorem.
\epf

Now we consider model spaces that are invariant under composition operators induced by linear fractional transformations. However, we will be restricted to inner functions that vanish at the origin. In fact, our main tool is the following adjoint formula \cite[Equation (3.6)]{Shapiro:Adjoint}, which concerns functions vanishing at the origin:

\blem\label{Shapiro:adjoint}
Let $\vp(z)=\frac{az+b}{cz+d}$ be a nonconstant holomorphic self-map of $\D$, where $ad - bc = 1$.
If $\sigma(z) = \frac{\bar{a}z-\bar{c}} {-\bar{b}z+\bar{d}}$, then
\[
C_\vp^*f(z)=z\sigma'(z) \frac{f(\sigma(z))}{\sigma(z)} \qquad (f\in H^2(\D), f(0)=0).
\]
\elem

As pointed our earlier, if $\theta$ is an inner function, then the constant function $1 \in \clq_\theta$ if and only if $\theta(0) = 0$. In the following, under the assumption that the inner functions vanish at the origin, we present a classification of model spaces that are invariant under composition operators induced by linear fractional transformations. Recall that $\clh(\D)$ is either $H^2(\D)$ or $H^\infty(\D)$.

\bthm\label{modelinv}
Let $\theta$ be an inner function, $\theta(0) = 0$, and let $\vp(z)= \frac{az+b}{cz+d}$ be a nonconstant holomorphic self-map of $\D$, where $ad - bc = 1$. Suppose $\sigma(z)= \frac{\bar{a}z-\bar{c}}{-\bar{b}z+\bar{d}}$ and
\[
\psi = \frac{\frac{\theta}{z}\circ\sigma}{\frac{\theta}{z}}.
\]
The following are equivalent:
\begin{enumerate}
  \item $\clq_\theta \in \text{Lat} C_\vp$.
  \item $\psi \in \clh(\D)$.
  \item $\frac{\theta}{z}H^2(\D) \in \text{Lat} C_\sigma$.
\end{enumerate}
\ethm
\bpf
Since $\theta(0)=0$, we see that $\frac{\theta}{z}$ is an inner function. Suppose that $\clq_\theta$ is invariant under $C_\vp$. Lemma \ref{eqn: Xf} implies that
\[
Xf = \frac{1}{\theta}C_\vp^*(\theta f) \qquad (f \in H^2(\D)),
\]
defines a bounded linear operator $X$ on $H^2(\D)$ and $C_{\vp}^* M_\theta = M_\theta X$. In view of Lemma \ref{Shapiro:adjoint}, it then follows that
\[
Xf = \sigma' \frac{z (\theta\circ\sigma)}{\theta\sigma} (f\circ\sigma) = \sigma' \psi (f\circ\sigma),
\]
where $\sigma'$ is the derivative of $\sigma$. In particular
\[
X1 = \sigma' \psi \in H^2(\D),
\]
from which, using the fact that $\frac{1}{\sigma'} \in H^\infty(\D)$, it follows that $\psi \in H^2(\D)$. And finally, by Theorem \ref{main}, $\ds\frac{\theta}{z}H^2(\D)$ is invariant under $C_\sigma$, which proves $(1) \Rightarrow (2) \Rightarrow (3)$.

\noindent Now if $(3)$ is true, then Theorem \ref{main} implies that $(2)$ is true, that is, $\psi \in H^\infty(\D)$. In particular, since $\sigma' \in H^\infty(\D)$, $X$ (defined as above) defines a bounded linear operator on $H^2(\D)$ and satisfy the equality $C_{\vp}^* M_\theta = M_\theta X$. By Lemma \ref{eqn: Xf}, $\clq_\theta$ is invariant under $C_\vp$, which complete the proof of $(2) \Rightarrow (1)$.
\epf

We do not know, in general, how to classify an inner function $\theta$ for which $\clq_\theta$ is invariant under a composition operator induced by a linear fractional transformation.

\section{Applications}\label{Mjay2-sec5}

The main goal of this section is to discuss some direct applications of our main results to some elementary invariant subspaces of composition operators.

First, we prove that the identity map is the only self-map of $\D$, which leaves all model spaces invariant.

\bthm
Let $\vp$ be a holomorphic self-map of $\D$. Then $\theta H^2(\D)$ is  invariant under $C_\vp$ for all inner function $\theta$ if and only if $\vp$ is the identity map.
\ethm
\bpf
If $\vp(z) = z$, $z \in \D$, then $C_\vp = I_{H^2(\D)}$ (the identity operator). In particular, $\theta H^2(\D)$ is invariant under $C_\vp$. Conversely, suppose $\theta H^2(\D)\in \mbox{Lat}C_\vp$ for all inner function $\theta$. Fix $\alpha \in \D$. Then $\alpha$ is the only zero of the Blaschke factor (which is also an inner function) $b_\alpha(z) = \frac{\alpha - z}{1-\bar{\alpha} z}$, $z \in \D$. Since $b_\alpha H^2(\D)\in \mbox{Lat}C_\vp$, by assumption, Theorem \ref{main} implies that
\[
\frac{b_\alpha\circ\vp}{b_\alpha}\in \es(\D).
\]
This forces $b_\alpha(\vp(\alpha))=0$, and hence $\vp(\alpha)=\alpha$. Since $\alpha \in \D$ is arbitrary, we have that $\vp(z)=z$, $z\in \mathbb{D}$, which completes the proof.
\epf

Before examining the other extreme property of the invariant subspaces of composition operators, we need to recall the following fact which concerns fixed points of holomorphic self-maps \cite[Corollary 3.4]{ISCO}: Let $\varphi$  be a holomorphic self-map of $\mathbb{D}$ and let $w \in \D$ be the fixed point of $\vp$. If $\theta$ is an inner function and $\theta(w)\neq 0$, then $\theta H^2(\D)$ is invariant under $C_\varphi$ if and only if $\theta\circ\varphi=\theta$.

\noindent Recall that a point $w \in \D$ is said to be a \textit{fixed point} of a holomorphic self-map $\vp$ of $\D$ if $\vp(w) = w.$

\bthm\label{thm: all phi}
Let $\theta$ be an inner function. Then $\theta H^2(\D) \in \text{Lat} C_{\vp}$ for all self-maps $\vp$ of $\D$ if and only if $\theta$ is an unimodular constant.
\ethm
\bpf
If $\theta$ is an unimodular constant, then $\frac{\theta\circ\vp}{\theta}$ is also unimodular constant. Theorem \ref{main} then implies that $\theta H^2(\D)\in \mbox{Lat} C_\vp$.

\noindent Conversely, assume that $\theta H^2(\D) \in \text{Lat} C_\vp$ for all holomorphic self-map $\vp$ of $\D$. First, we claim that $\theta$ has no zeros in $\D$. Suppose towards a contradiction that $\theta(\alpha)=0$ for some $\alpha \in \D$. Since $\theta$ is an inner function, there exists $\beta\in \D$ such that  $\theta(\beta) \neq 0$. Set
\[
\vp = b_{\beta} \circ b_{\alpha}.
\]
Note that $\vp$ is an inner function and $\vp(\alpha) = \beta$. Therefore, $\frac{\theta\circ\vp}{\theta} \notin H^\infty(\D)$, and hence, by Theorem \ref{main}, $\theta H^2(\D)$ is not invariant under $C_\vp$. This contradiction completes the proof of the claim that $\theta(z)\neq 0$ for all $z\in \D$. If we consider a constant function $\vp\equiv c$, then in view of the discussion preceding the statement of this theorem, $C_\vp (\theta H^2(\D)) \subseteq \theta H^2(\D)$ implies that
\[
\theta=\theta\circ\vp\equiv\theta(c).
\]
Since $\theta$ is an inner function, $\theta(c)$ must be unimodular, and hence $\theta$ is an unimodular constant. This completes the proof of the theorem.
\epf

Now we turn to special model invariant subspaces of composition operators. In the following, we will let $\theta$ and $\vp$ denote inner function and holomorphic self-map of $\D$, respectively.

\bthm\label{thm: all quotient inv sub}
The following holds true:
\begin{enumerate}
\item $\clq_\theta \in \text{Lat} C_\vp$ for all inner function $\theta$ if and only if $\vp$ is the identity map.
\item $\clq_\theta \in \text{Lat} C_\vp$ for all self-map $\vp$ if and only if $\theta (z)=\alpha z$ or $\theta\equiv\alpha$, where $\alpha$ is an unimodular constant.
\end{enumerate}
\ethm
\bpf
If $\clq_\theta \in \text{Lat} C_\vp$ for all inner function $\theta$, then, in particular, $\clq_{b_{\alpha}}, \clq_{b_{\beta}} \in \text{Lat} C_\vp$, where $\alpha \neq \beta$ are in $\D \setminus \{0\}$. By Theorem \ref{affine}, there exist nonzero scalars $a$ and $b$ such that
\[
\vp(z)= \frac{1-a}{\bar{\alpha}}+ a z = \frac{1-b}{\bar{\beta}} + b z \qquad (z \in \D).
\]
Thus, comparing both sides of the above and noting the fact that $\alpha \neq \beta$, we have $a = b =1$, which implies that $\vp(z)=z$, $z\in \mathbb{D}$. The converse part of $(1)$ is trivial.

\noindent Now we turn to prove $(2)$. If $\theta (z)=\alpha z$ for some unimodular $\alpha$, then $\clq_\theta = k(\cdot, 0) \mathbb{C}$ is the space of all constant functions, which is clearly invariant under all composition operators (cf. Remark \ref{constant}). In the case of $\theta\equiv\alpha$, we have $\clq_\theta=\{0\}$, which is also invariant under all composition operators.

\noindent For the converse part, we first note that if $\clq_\theta=\{0\}$, then $\theta H^2(\D)=H^2(\D)$, and hence $\theta$ is an unimodular constant. So suppose $\clq_\theta\neq \{0\}$. Then there exist $\alpha \in \D$ and $f\in \clq_\theta$ such that $f(\alpha)\neq 0$. By taking $\vp\equiv \alpha$, we see that $1 \in \clq_\theta$, which implies
\[
\theta(0)=\langle \theta, 1\rangle=0.
\]
Then $\omega = \frac{\theta}{z}$ is an inner function. By Theorem \ref{modelinv}, $\omega H^2(\D)$ invariant under $C_\sigma$ for every linear fractional self-map $\sigma$ of $\D$.

\noindent Now we claim that $\omega$ has no zeros in $\D$. Suppose towards a contradiction that $\omega(\alpha)=0$ for some $\alpha \in \D$. Since $\omega$ is an inner function, $\omega(\beta) \neq 0$ for some $\beta \in \D$. Set $\sigma = b_\beta \circ b_\alpha$ so that $\sigma(\alpha) = \beta$. Then
\[
\frac{\omega\circ\sigma}{\omega} \notin H^\infty(\D),
\]
and hence $\omega H^2(\D)$ is not invariant under $C_\sigma$. This is a contradiction, and hence $\omega(z)\neq 0$, $z\in \D$.

\noindent Finally, consider a non-automorphic and holomorphic self-map $\sigma$ of $\D$ with a fixed point in $\D$ (for instance, consider $\sigma(z)= \frac{1}{2}z$). As $\omega H^2(\D)$ invariant under $C_\sigma$, by \cite[Corollary 3.5]{ISCO}, it follows that $\omega$ is an unimodular constant. This implies that $\theta (z)=\alpha z$ for some unimodular constant $\alpha$.
\epf

We conclude this section by proving a reducing property of composition operators. Here we will be restricted to finite-dimensional model spaces, as in Theorem \ref{affine}. Before we begin this discussion, we recall the following result \cite[Corollary 2.4]{ISCO}:

\bthm  \label{Blaschke}
Let $B$ be a Blaschke product and let $\vp$ be a holomorphic self-map of $\mathbb{D}$. Then $BH^2(\D) \in \text{Lat}C_\vp$ if and only if $\text{multiplicity}_w B \leq \text{multiplicity}_w {B\circ\vp}$ for all $w \in \D$ such that $B(w)=0$.
\ethm

We are now ready for the reducing subspace property of elementary composition operators.

\bthm\label{thm: reducing}
Let $\vp$ be a holomorphic self-map of $\D$, $\alpha \in \D$, $n \geq 1$, and suppose $\theta(z)=\left(\frac{z- \alpha}{1-\bar{\alpha}z}\right)^n$. Then $\clq_\theta$ reduces $C_\vp$ if and only if there exist a scalar $c$, $|c| \leq 1$, and $\psi \in \cls(\D)$ such that
\[
\vp(z)= \begin{cases}
z \psi & \mbox{if } \alpha = 0 \text{ and } n = 1 \\
c z & \mbox{if } \alpha = 0 \text{ and } n \geq 2\\
z & \mbox{if } \alpha \neq 0.
\end{cases}
\]
\ethm
\bpf
Suppose $\alpha = 0$. If $n=1$, then the one-dimensional subspace $(z H^2(\D))^{\bot} = \clq_z$ is invariant under every composition operator (see Remark \ref{constant}), where $z H^2(\D)$  is invariant under $C_\vp$ if and only if $\vp= z \psi$ for some $\psi \in \es(\D)$ (see Theorem \ref{Blaschke}). If $n \geq 2$, then, by Theorem \ref{affine}, $\clq_\theta$ is invariant under $C_\vp$ if and only if $\vp$ is an affine map. By Theorem \ref{Blaschke}, $\theta H^2(\D)$ is invariant under $C_\vp$ if and only if $\vp= z \psi$ for some $\psi \in \es(\D)$. In particular, $\vp(0)=0$, and hence $\theta H^2(\D)$ reduces $C_\vp$ if and only if $\vp(z)= cz$ for some $|c|\leq 1$.

\noindent Finally, assume that $\alpha \neq 0$. If $\vp$ is the identity map, then $C_\vp = I_{H^2(\D)}$, and hence $\clq_\theta$ reduces $C_\vp$. Next, we assume that $\clq_\theta$ reduces $C_\vp$. By Theorem \ref{affine}, there exists a scalar $c (\neq 0)$ such that $\vp(z)= \frac{1-c}{\bar{\alpha}} + cz$. Again, by Theorem \ref{Blaschke}, $(\theta\circ\vp)(\alpha)=0$, and hence $\vp(\alpha)= \alpha$. Now
\[
\frac{1-c}{\bar{\alpha}} + c \alpha = \alpha,
\]
is equivalent to
\[
(1-c)(1-|\alpha|^2)=0.
\]
As $|\alpha| < 1$, we must have $c=1$, that is, $\vp(z)=z$, $z\in \D$, which completes the proof.
\epf

Finally, we would like to point out that the present theory of model spaces and composition operators complement the recent paper \cite{ISCO}. The classification of Beurling type invariant subspaces of composition operators (\cite[Theorem 2.3]{ISCO} or Theorem \ref{main}) also plays an important role in our consideration, and evidently, some of our results are in duality (cf. Theorem \ref{tm: FLT 1} and \ref{modelinv}). On the other hand, model spaces are fairly complex objects in operator theory and function theory, which also represent a large class of bounded linear operators. The present analysis clearly shows a challenging role of model spaces to the theory of composition operators (and vice versa).

\smallskip

\subsection*{Acknowledgments}
The first author thanks the National Board for Higher Mathematics (NBHM), India, for providing financial support  to carry out this research. The second author is supported in part by NBHM (National Board of Higher Mathematics, India) grant NBHM/R.P.64/2014, and the Mathematical Research
Impact Centric Support (MATRICS) grant, File No : MTR/2017/000522, by the Science and Engineering Research Board (SERB), Department of Science \& Technology (DST), Government of India.

\end{document}